\documentclass[12pt]{article}
\usepackage{amsmath}
\usepackage{amssymb}
\usepackage{vatola}

\def\q{\quad}
\def\qq{\qquad}
\def\mod{\pmod}
\def\t{\text}
\def\f{\frac}
\def\e{\equiv}

\def\phq#1{\varphi(q^{#1})}
\def\psq#1{\psi(q^{#1})}
\def\qtq#1{\q\t{#1}\q}

\let \pro=\proclaim
\let \endpro=\endproclaim

\begin{document}
 \par\q\newline
\centerline {\bf  Ramanujan's theta functions and linear
combinations of three triangular numbers}

$$\q$$
\centerline{Zhi-Hong Sun}
\par\q\newline
\centerline{School of Mathematical Sciences}
 \centerline{Huaiyin
Normal University} \centerline{Huaian, Jiangsu 223300, P.R. China}
\centerline{Email: zhsun@hytc.edu.cn} \centerline{Homepage:
http://www.hytc.edu.cn/xsjl/szh}

 \abstract{Let $\Bbb Z$ be the set of integers. For positive
 integers $a,b,c$ and $n$ let $N(a,b,c;n)$ be the number of
 representations of $n$ by $ax^2+by^2+cz^2$, and
 let $t(a,b,c;n)$ be the number of
 representations of $n$ by
 $ax(x+1)/2+by(y+1)/2+cz(z+1)/2
 $ $(x,y,z\in\Bbb Z)$. In this paper, by using
 Ramanujan's theta functions $\varphi(q)$ and $\psi(q)$
 we reveal the relation between $t(2,3,3;n)$ and $N(1,3,3;n+1)$, and
 the relation between $t(1,1,6;n)$ and $N(1,1,3;n+1)$.
We also obtain formulas for $t(a,3a,4b;n),$
$t(a,7a,4b;n),t(3a,5a,4b;n)$ and $t(a,15a,4b;n)$ under certain
congruence conditions, where $a$ and $b$ are positive odd integers.
  In addition, we pose many conjectures on $t(a,b,c;n)$ for some special values of
$(a,b,c)$.

 \par\q
 }
 \endabstract

\section*{1. Introduction}

\par  Let $\Bbb Z$, $\Bbb Z^+$ and $\Bbb N$ be the set of
integers, the set of positive integers
 and the set of nonnegative integers, respectively,
 and let
 $\Bbb Z^3=\Bbb Z\times
  \Bbb Z\times\Bbb Z$ and $\Bbb N^3=\Bbb N\times
  \Bbb N\times\Bbb N$.
  For $a,b,c\in\Bbb Z^+$ and $n\in\Bbb N$  set
$$\align &N(a,b,c;n)=\big|\{(x,y,z)\in \Bbb Z^3\ |
\ n=ax^2+by^2+cz^2 \}\big|,
\\&t(a,b,c;n)=\Big|\Big\{(x,y,z)\in \Bbb
Z^3\ \big|\ n\ =a\f{x(x+1)}2+ b\f{y(y+1)}2+c\f{z(z+1)}2\Big\}\Big|,
\\&T(a,b,c;n)=\Big|\Big\{(x,y,z)\in \Bbb
N^3\ \big|\ n\ =a\f{x(x+1)}2+ b\f{y(y+1)}2+c\f{z(z+1)}2\Big\}\Big|.
\endalign$$
Since $\f{x(x+1)}2=\f{(-1-x)(-x)}2$ we see that
$$t(a,b,c;n)=8T(a,b,c;n).\tag 1.1$$
For $a,b,c\in\Bbb Z^+$ let
$$C(a,b,c)=\f{i_1(i_1-1)(i_1-2)(i_1-3)}4
+\f{i_1(i_1-1)i_2}2+i_1i_3,\tag 1.2$$ where $i_j$ denotes the number
of elements in $\{a,b,c\}$ which are equal to $j$.
 In 2005 Adiga, Cooper and
Han [ACH] showed that
$$t(a,b,c;n)=\f{2}{2+C(a,b,c)}
N(a,b,c;8n+a+b+c) \q\t{for $a+b+c\le 7$}. \tag 1.3$$ In 2008 Baruah,
Cooper and Hirschhorn [BCH] proved that
$$\aligned &t(a,b,c;n)
\\&=\f{2}{2+C(a,b,c)}(N(a,b,c;8n+8)-N(a,b,c;2n+2))
\q \t{for $a+b+c=8$}.\endaligned \tag 1.4$$

\par Ramanujan's theta functions $\varphi(q)$ and $\psi(q)$ are defined
by
$$\varphi(q)=\sum_{n=-\infty}^{\infty}q^{n^2}=1+2\sum_{n=1}^{\infty}
q^{n^2}\qtq{and} \psi(q)=\sum_{n=0}^{\infty}q^{n(n+1)/2}\
(|q|<1).\tag 1.5$$ It is evident that for positive integers
$a_1,\ldots,a_k$ and $|q|<1$,
$$\align&\sum_{n=0}^{\infty}N(a_1,\ldots,a_k;n)q^{n}=\varphi(q^{a_1})
\cdots\varphi(q^{a_k}),\tag 1.6
\\&\sum_{n=0}^{\infty}t(a_1,\ldots,a_k;n)q^{n}=2^k\psi(q^{a_1})
\cdots \psi(q^{a_k}).\tag 1.7\endalign$$ There are many identities
involving $\varphi(q)$ and $\psi(q)$.  From [BCH, Lemma 4.1] or [Be]
we know that for $|q|<1$,
 $$\align &\psi(q)^2=\varphi(q)\psi(q^2),\tag
 1.8\\&\varphi(q)=\varphi(q^4)+2q\psi(q^8)=\varphi(q^{16})+2q^{4}\psi(q^{32})
  +2q\psi(q^{8}),\tag 1.9
 \\&\varphi(q)^2=\phq 2^2+4q\psq 4^2=\phq 4^2+4q^2\psq 8^2+4q\psq 4^2,\tag 1.10
  \\&\psi(q)\psi(q^3)=\varphi(q^6)\psi(q^4)+q
 \varphi(q^2)\psi(q^{12}).\tag 1.11
 \endalign$$
By [S1, Lemma 2.4],
 $$ \varphi(q)^2=\phq 8^2+4q^4\psq{16}^2+4q^2\psq 8^2+4q\phq{16}\psq 8
 +8q^5\psq 8\psq{32}.
\tag 1.12$$ By [S1, Lemma 2.3], for $|q|<1$ we have
$$\aligned\varphi(q)\phq 3&=
\phq {16}\phq{48}+4q^{16}\psq{32}\psq{96}+2q \phq {48}\psq
8+2q^3\phq{16}\psq{24} \\&\q+6q^4\psq 8\psq{24}+4q^{13}\psq
8\psq{96}+4q^7\psq{24}\psq{32}.\endaligned \tag 1.13$$

 \par Using theta function identities we may establish some relations
 between $t(a,b,c;n)$ and $N(a,b,c;8n+a+b+c)$ for some certain values of
 $(a,b,c)$. See [S2,S3].
 \par In this paper we reveal
the relation between $t(2,3,3;n)$ and $N(1,3,3;n+1)$, and
 the relation between $t(1,1,6;n)$ and $N(1,1,3;n+1)$.
We also obtain formulas for $t(a,3a,4b;n),$
$t(a,7a,4b;n),t(3a,5a,4b;n)$ and $t(a,15a,4b;n)$ under certain
congruence conditions, where $a$ and $b$ are positive odd integers.
  In addition, we pose
many conjectures on $t(a,b,c;n)$ for some special values of
$(a,b,c)$.

\section*{2. The relation between $t(2,3,3;n)$ and $N(1,3,3;n+1)$}
\par By (1.4), for  $n\in\Bbb Z^+$ we have
$$t(2,3,3;n)=N(2,3,3;8n+8)-N(2,3,3;2n+2).$$
By [S2, Theorem 6.1], if $a,b,c,n\in\Bbb Z^+$, $2\nmid ab$, $4\mid
a-b$ and $4\mid c-2$, then
$$t(a,b,c;n)=N(a,b,c;8n+a+b+c)-N(a,b,c;2n+(a+b+c)/4).$$
\par Now we present the following formula for $t(2,3,3;n)$.

\pro{Theorem 2.1} For $m\in\Bbb N$ we have
$$\align &N(2,3,3;2m)=N(1,3,3;m),
\\&t(2,3,3;8m+6)=2t(1,3,3;m),\q t(2,3,3;32m+27)=4t(1,3,3;m).
\endalign$$
For $n\in\Bbb Z^+$ with $n\not\e 15\mod{16}$ we have
$$t(2,3,3;n)=\cases  4N(1,3,3;n+1)&\t{if $n\e 0,1\mod 4$,}
\\ 2N(1,3,3;n+1)&\t{if $n\e 2\mod 8$,}
\\ N(1,3,3;n+1)&\t{if $n\e 6\mod 8$ or $n\e 27\mod {32}$,}
\\ \f 85N(1,3,3;n+1)&\t{if $n\e 3,7\mod {16}$,}
\\ \f 43N(1,3,3;n+1)&\t{if $n\e 11\mod {32}$.}
\endcases$$
\endpro
Proof. By (1.6) and (1.9),
$$\sum_{n=0}^{\infty}N(2,3,3;n)q^n
=\phq 2\phq 3^2=\phq 2(\phq{6}^2+4q^3\psq{12}^2).$$ Thus,
$$\sum_{m=0}^{\infty}N(2,3,3;2m)q^{2m}=\phq 2\phq 6^2$$
and so
$$\sum_{m=0}^{\infty}N(2,3,3;2m)q^m
=\varphi(q)\phq{3}^2=\sum_{m=0}^{\infty}N(1,3,3;m)q^m.$$ This yields
$N(2,3,3;2m)=N(1,3,3;m)$.
 By (1.12),
$$\aligned \varphi(q^3)^2&
=\phq{24}^2+4q^{12}\psq{48}^2 +4q^3\phq{48}\psq
{24}\\&\qq+8q^{15}\psq {24}\psq{96} +4q^6\psq {24}^2.\endaligned\tag
2.1$$ Thus, using (1.9), (1.10) and the above we see that
$$\aligned&\sum_{m=0}^{\infty}N(1,3,3;m)q^m\\&=\varphi(q)\phq 3^2
=(\phq 4+2q\psq 8)(\phq{12}^2+4q^6\psq{24}^2+4q^3\psq{12}^2) \\&
=(\phq {16}+2q\psq 8+2q^4\psq{32}) \big(\phq{24}^2+4q^{12}\psq{48}^2
\\&\qq\qq\qq\qq+4q^3\phq{48}\psq {24}+8q^{15}\psq {24}\psq{96}
+4q^6\psq {24}^2\big).\endaligned\tag 2.2$$ This yields
$$\align&\sum_{m=0}^{\infty}N(1,3,3;4m+1)q^{4m+1}=2q\psq
8\phq{12}^2,
\\&\sum_{m=0}^{\infty}N(1,3,3;4m+2)q^{4m+2}=4q^6\psq{24}^2\phq 4
\endalign$$ and so
$$\align &\sum_{m=0}^{\infty}N(1,3,3;4m+1)q^m=2\phq 3^2\psq 2,\tag
2.3
\\&\sum_{m=0}^{\infty}N(1,3,3;4m+2)q^m=4q\varphi(q)\psq 6^2.\tag 2.4
\endalign$$
Also, from (2.2) we deduce that
$$\align &\sum_{m=0}^{\infty}N(1,3,3;8m+3)q^{8m+3}
=\phq{16}\cdot 4q^3\phq{48}\psq{24}+2q^4\psq{32}\cdot
8q^{15}\psq{24}\psq{96},
\\&\sum_{m=0}^{\infty}N(1,3,3;8m+7)q^{8m+7}
\\&\qq=\phq{16}\cdot 8q^{15}\psq{24}\psq{96}+2q\psq 8\cdot 4q^6\psq{24}^2
+2q^4\psq{32}\cdot 4q^3\phq{48}\psq{24}
\endalign$$
and so
$$\align &\sum_{m=0}^{\infty}N(1,3,3;8m+3)q^m
=4\phq{2}\phq{6}\psq{3}+16q^2\psq{4}\cdot \psq{3}\psq{12},\tag 2.5
\\&\sum_{m=0}^{\infty}N(1,3,3;8m+7)q^m\tag 2.6
\\&\q=8q\phq{2}\psq{3}\psq{12}+8\psi(q)\psq{3}^2
+8\psq{4}\phq{6}\psq{3}=16\psi(q)\psq 3^2.
\endalign$$
On the other hand,
$$\aligned \sum_{n=0}^{\infty}t(2,3,3;n)q^n
&=8\psq 2\psq 3^2=8\psq 2\phq 3\psq 6
\\&=8(\phq{12}+2q^3\psq{24})(\phq{12}\psq 8+q^2\phq 4\psq{24}).
\endaligned\tag 2.7$$
Hence
$$\align &\sum_{m=0}^{\infty}t(2,3,3;4m)q^{4m}=8\phq{12}^2\psq 8,
\\&\sum_{m=0}^{\infty}t(2,3,3;4m+1)q^{4m+1}=8\cdot 2q^3\psq{24}\cdot
q^2\phq 4\psq{24},
\\&\sum_{m=0}^{\infty}t(2,3,3;4m+2)q^{4m+2}
=8\phq{12}\cdot q^2\phq 4\psq{24},
\\&\sum_{m=0}^{\infty}t(2,3,3;4m+3)q^{4m+3}=16q^3\psq{24}\phq{12}\psq
8
\endalign$$ and so
$$\align &\sum_{m=0}^{\infty}t(2,3,3;4m)q^m=8\phq{3}^2\psq 2,\tag
2.8
\\&\sum_{m=0}^{\infty}t(2,3,3;4m+1)q^m=16q\varphi(q)\psq{6}^2,\tag
2.9
\\&\sum_{m=0}^{\infty}t(2,3,3;4m+2)q^m=8\varphi(q)\phq 3\psq 6,\tag
2.10
\\&\sum_{m=0}^{\infty}t(2,3,3;4m+3)q^m=16\phq 3\psq 2\psq 6.\tag
2.11
\endalign$$ Combining (2.3),(2.4), (2.7) and (2.8) yields
$$t(2,3,3;4m)=4N(1,3,3;4m+1),\q t(2,3,3;4m+1)=4N(1,3,3;4m+2).$$
Hence $t(2,3,3;n)=4N(1,3,3;n+1)$ for $n\e 0,1\mod 4$.
\par Now we consider the case $n\e 2,6\mod 8$. By (2.10),
$$\sum_{m=0}^{\infty}t(2,3,3;4m+2)q^m=
8(\phq 4+2q\psq 8)(\phq{12}+2q^3\psq{24})\psq 6.$$ Thus,
$$\align &\sum_{m=0}^{\infty}t(2,3,3;8m+2)q^{2m}
=8(\phq 4\phq {12}+4q^4\psq 8\psq{24})\psq 6,
\\&\sum_{m=0}^{\infty}t(2,3,3;8m+6)q^{2m+1}
=8\psq 6(2q\psq 8\phq{12}+\phq 4\cdot 2q^3\psq{24})
\endalign$$ and so
$$\align&\sum_{m=0}^{\infty}t(2,3,3;8m+2)q^m
=8(\phq 2\phq {6}+4q^2\psq 4\psq{12})\psq 3,\tag 2.12
\\&\sum_{m=0}^{\infty}t(2,3,3;8m+6)q^m
=16\psq 3(\psq 4\phq{6}+q\phq 2\psq{12})\tag 2.13
\\&\qq\qq\qq\qq\qq\;=16\psi(q)\psq
3^2=2\sum_{m=0}^{\infty}t(1,3,3;m)q^m.
\endalign$$
Thus,
$$t(2,3,3;8m+6)=2t(1,3,3;m).$$
 Combining (2.5),(2.6),(2.12) and (2.13) yields
$$t(2,3,3;8m+2)=2N(1,3,3;8m+3),\q t(2,3,3;8m+6)=N(1,3,3;8m+7).$$
That is,
$$t(2,3,3;n)=\cases 2N(1,3,3;n+1)&\t{if $n\e 2\mod 8$,}
\\N(1,3,3;n+1)&\t{if $n\e 6\mod 8$.}\endcases$$
\par By (2.11), $$\align &\sum_{m=0}^{\infty}t(2,3,3;4m+3)q^m
\\&=16\phq 3\psq 2\psq 6=16(\phq{12}+2q^3\psq{24})(\phq{12}\psq
8+q^2\phq 4\psq{24})
\\&=16(\phq{12}^2\psq 8+2q^5\phq 4\psq{24}^2+q^2\phq
4\phq{12}\psq{24}+2q^3\phq{12}\psq 8\psq{24}).\endalign$$ Hence
$$\align&\sum_{m=0}^{\infty}t(2,3,3;16m+3)q^{4m}
=16\phq{12}^2\psq 8,
\\&\sum_{m=0}^{\infty}t(2,3,3;16m+7)q^{4m+1}
=32q^5\phq 4\psq {24}^2,
\\&\sum_{m=0}^{\infty}t(2,3,3;16m+11)q^{4m+2}
=16q^2\phq 4\phq{12}\psq{24}.
\endalign$$ Therefore,
$$\align&\sum_{m=0}^{\infty}t(2,3,3;16m+3)q^m
=16\phq{3}^2\psq 2,\tag 2.14
\\&\sum_{m=0}^{\infty}t(2,3,3;16m+7)q^m
=32q\varphi(q)\psq {6}^2,\tag 2.15
\\&\sum_{m=0}^{\infty}t(2,3,3;16m+11)q^m
=16\varphi(q)\phq 3\psq 6 .\tag 2.16\endalign$$ From (2.2) we see
that
$$\sum_{m=0}^{\infty}N(1,3,3;4m)q^{4m}=\phq 4\phq{12}^2+8q^4\psq
8\psq{12}^2$$ and so
$$\align &\sum_{m=0}^{\infty}N(1,3,3;4m)q^{m}
\\&=\varphi(q)\phq{3}^2+8q\psq 2\psq{3}^2
\\&=\varphi(q)\phq{3}^2+8q\phq 3\psq 2\psq 6
\\&=\varphi(q)\phq{3}^2+8q\phq 3(\phq{12}\psq 8+q^2\phq 4\psq{24})
\\&=(\phq 4+2q\psq 8)(\phq{12}+2q^3\psq{24})^2
\\&\q+8q(\phq{12}+2q^3\psq{24}) (\phq{12}\psq 8+q^2\phq 4\psq{24})
\\&=\phq 4\phq{12}^2+24q^4\phq{12}\psq 8\psq{24}
+10q\phq{12}^2\psq 8
\\&\q+20q^6\phq 4\psq{24}^2+4q^3(3\phq 4\phq{12}\psq{24}
+2q^4\psq 8\psq{24}^2).
\endalign$$
Therefore,
$$\align
&\sum_{m=0}^{\infty}N(1,3,3;4(4m+1))q^{4m+1}=10q\phq{12}^2\psq 8,
\\&\sum_{m=0}^{\infty}N(1,3,3;4(4m+2))q^{4m+2}=20q^6\phq 4\psq
{24}^2,
\\&\sum_{m=0}^{\infty}N(1,3,3;4(4m+3))q^{4m+3}
=4q^3(3\phq 4\phq{12}\psq{24} +2q^4\psq 8\psq{24}^2)\
\endalign$$
and so
$$\align
&\sum_{m=0}^{\infty}N(1,3,3;16m+4)q^m=10\phq{3}^2\psq 2,\tag 2.17
\\&\sum_{m=0}^{\infty}N(1,3,3;16m+8)q^m=20q\varphi(q)\psq
{6}^2,\tag 2.18
\\&\sum_{m=0}^{\infty}N(1,3,3;16m+12)q^m =4(3\varphi(q)
\phq{3}\psq{6} +2q\psq 2\psq{6}^2).\tag 2.19
\endalign$$
Combining (2.17)-(2.18) with (2.14) and (2.15) yields
$$t(2,3,3;16m+3)=\f 85N(1,3,3;16m+4),\q
t(2,3,3;16m+7)=\f 85N(1,3,3;16m+8).$$ That is,
$$t(2,3,3;n)=\f 85N(1,3,3;n+1)\qtq{for}n\e 3,7\mod{16}.$$
From (2.16) we see that
$$\align &\sum_{m=0}^{\infty}t(2,3,3;16m+11)q^m
\\&=16\varphi(q)\phq 3\psq 6
=16(\phq 4+2q\psq 8)(\phq{12}+2q^3\psq{24})\psq 6
\\&=16\psq 6(\phq 4\phq{12}+4q^4\psq 8\psq{24}+2q(\phq{12}\psq 8
+q^2\phq 4\psq{24}))
\\&=16\psq 6(\phq 4\phq{12}+4q^4\psq 8\psq{24})+32q\psq 2\psq 6^2.
\endalign$$
Thus,
$$\align&\sum_{m=0}^{\infty}t(2,3,3;32m+11)q^{2m}
=16\psq 6(\phq 4\phq{12}+4q^4\psq 8\psq{24}),
\\&\sum_{m=0}^{\infty}t(2,3,3;32m+27)q^{2m+1}
=32q\psq 2\psq 6^2.\endalign$$ Hence
$$\align&\sum_{m=0}^{\infty}t(2,3,3;32m+11)q^m
=16\psq 3(\phq 2\phq{6}+4q^2\psq 4\psq{12}),\tag 2.20
\\&\sum_{m=0}^{\infty}t(2,3,3;32m+27)q^m
=32\psi(q)\psq 3^2=4\sum_{m=0}^{\infty}t(1,3,3;m)q^m.\tag
2.21\endalign$$ Thus, $t(2,3,3;32m+27)=4t(1,3,3;m).$ From (2.19) we
see that
$$\align&\sum_{m=0}^{\infty}N(2,3,3;16m+12)q^m
\\&=4(3\varphi(q) \phq{3}\psq{6} +2q\psq 2\psq 6^2)
\\&=12\psq 6(\phq 4+2q\psq 8)(\phq{12}+2q^3\psq{24})+8q\psq 2\psq 6^2
.\endalign$$ Thus,
$$\align &\sum_{m=0}^{\infty}N(1,3,3;32m+12)q^{2m}
\\&\q=12\psq 6(\phq 4\phq{12}+4q^4\psq 8\psq{24}),
\\&\sum_{m=0}^{\infty}N(2,3,3;2(32m+28))q^{2m+1}
\\&\q=12\psq 6(2q\psq 8\phq{12}+2q^3\phq 4\psq{24})+8q\psq 2\psq 6^2
\\&\q=12\psq 6\cdot 2q\psq 2\psq 6+8q\psq 2\psq 6^2=32q\psq 2\psq 6^2.
\endalign$$ Therefore,
$$\align &\sum_{m=0}^{\infty}N(1,3,3;32m+12)q^m
=12\psq 3(\phq 2\phq{6}+4q^2\psq 4\psq{12}),
\\&\sum_{m=0}^{\infty}N(1,3,3;32m+28)q^m
=32\psi(q)\psq 3^2.\endalign$$ This together with (2.20) and (2.21)
gives
$$\align&
t(2,3,3;32m+27)=N(1,3,3;32m+28),\\& t(2,3,3;32m+11)=\f
43N(1,3,3;32m+12).\endalign$$
 This completes the proof.

 \section*{3. The relation between $t(1,1,6;n)$ and $N(1,1,3;n+1)$}
 \par By (1.4), for $n\in\Bbb Z^+$ we have
 $$t(1,1,6;n)=N(1,1,6;8n+8)-N(1,1,6;2n+2).$$
\pro{Theorem 3.1} For  $n\in\Bbb Z^+$ we have
$$\align &N(1,1,6;2n)=N(1,1,3;n),\q t(1,1,6;8n+4)=2t(1,1,3;n),
\\&t(1,1,6;32n+19)=4t(1,1,3;n).\endalign$$
For $n\in\Bbb Z^+$ with $n\not\e 15\mod {16}$ we have
$$t(1,1,6;n)=\cases
 4N(1,1,3;n+1)&\t{if $n\e 1,2\mod 4$,}
\\2N(1,1,3;n+1)&\t{if $n\e 0\mod 8$,}
\\N(1,1,3;n+1)&\t{if $n\e 4\mod 8$ or $n\e 19\mod{32}$,}
\\\f 85N(1,1,3;n+1)&\t{if $n\e 7,11\mod {16}$,}
\\\f 43N(1,1,3;n+1)&\t{if $n\e 3\mod {32}$.}
\endcases$$
\endpro
Proof. By (1.7)-(1.11),
$$\align &\sum_{n=0}^{\infty}t(1,1,6;n)q^n
\\&=8\psi(q)^2\psq 6=16\varphi(q)\psq 2\psq 6
\\&=8(\phq 4+2q\psq 8)(\phq{12}\psq 8+q^2\phq 4\psq{24})
\\&=8(\phq 4\phq{12}\psq 8+2q\phq{12}\psq 8^2+q^2\phq 4^2\psq{24}
+2q^3\phq 4\psq 8\psq{24}).
\endalign$$
Thus,
$$\align &\sum_{n=0}^{\infty}t(1,1,6;4n)q^{4n}=8\phq 4\phq{12}\psq 8,
\\ &\sum_{n=0}^{\infty}t(1,1,6;4n+1)q^{4n+1}=16q\phq{12}\psq 8^2,
\\&\sum_{n=0}^{\infty}t(1,1,6;4n+2)q^{4n+2}=8q^2\phq 4^2\psq{24},
\\&\sum_{n=0}^{\infty}t(1,1,6;4n+3)q^{4n+3}=16q^3\phq 4\psq 8\psq{24}).
\endalign$$
This yields
$$\align
&\sum_{n=0}^{\infty}t(1,1,6;4n)q^{n}=8\varphi(q)\phq{3}\psq 2,\tag
3.1
\\ &\sum_{n=0}^{\infty}t(1,1,6;4n+1)q^n=16\phq{3}\psq 2^2,\tag 3.2
\\&\sum_{n=0}^{\infty}t(1,1,6;4n+2)q^n=8\varphi(q)^2\psq{6},\tag
3.3
\\&\sum_{n=0}^{\infty}t(1,1,6;4n+3)q^n=16\varphi(q)\psq
2\psq{6}.\tag 3.4
\endalign$$
Since $\varphi(q)=\phq 4+2q\psq 8$, from (3.1) we see that
$$\sum_{n=0}^{\infty}t(1,1,6;4n)q^n=8\psq 2(\phq 4+2q\psq
8)(\phq{12}+2q^3\psq{24}).$$ This yields
$$\sum_{n=0}^{\infty}t(1,1,6;8n)q^{2n}=8\psq 2(\phq 4\phq{12}
+4q^4\psq 8\psq{24})$$ and
$$\sum_{n=0}^{\infty}t(1,1,6;8n+4)q^{2n+1}=16\psq 2(\phq{12}\psq
8+q^2\phq 4\psq{24})=32q\psq 2^2\psq 6.$$
Therefore,
$$\sum_{n=0}^{\infty}t(1,1,6;8n)q^n=8\psi(q)(\phq 2\phq{6}
+4q^2\psq 4\psq{12})\tag 3.5$$ and
$$\sum_{n=0}^{\infty}t(1,1,6;8n+4)q^n=16\psi(q)^2\psq 3.\tag 3.6$$
By (1.7) and (3.6), $t(1,1,6;8n+4)=2t(1,1,3;n)$.
 By (3.4),
$$\align&\sum_{n=0}^{\infty}t(1,1,6;4n+3)q^n
\\&=16\varphi(q)\psq
2\psq{6}=16(\phq 4+2q\psq 8)(\phq{12}\psq 8+q^2\phq 4\psq{24})
\\&=16(\phq 4\phq{12}\psq 8+2q\phq{12}\psq 8^2+q^2\phq 4^2\psq{24}
+2q^3\phq 4\psq 8\psq{24}).
\endalign$$
Thus
$$\align&\sum_{n=0}^{\infty}t(1,1,6;16n+3)q^{4n}
=16\phq 4\phq{12}\psq 8,
\\&\sum_{n=0}^{\infty}t(1,1,6;4(4n+1)+3)q^{4n+1}
=32q\phq{12}\psq 8^2,
\\&\sum_{n=0}^{\infty}t(1,1,6;4(4n+2)+3)q^{4n+2}
=16q^2\phq 4^2\psq{24}.\endalign$$ Replacing $q$ with $q^{1/4}$
yields
$$\align
&\sum_{n=0}^{\infty}t(1,1,6;16n+3)q^{n}=16\varphi(q)\phq 3\psq
2,\tag 3.7
\\&\sum_{n=0}^{\infty}t(1,1,6;16n+7)q^n
=32\phq{3}\psq 2^2,\tag 3.8
\\&\sum_{n=0}^{\infty}t(1,1,6;16n+11)q^n
=16\varphi(q)^2\psq{6}.\tag 3.9\endalign$$ By (3.7),
$$\sum_{n=0}^{\infty}t(1,1,6;16n+3)q^n
=16(\phq 4+2q\psq 8)(\phq {12}+2q^3\psq{24})\psq 2.$$ Thus,
$$\sum_{n=0}^{\infty}t(1,1,6;32n+3)q^{2n}
=16(\phq 4\phq{12}+4q^4\psq 8\psq{24})\psq 2$$ and
$$\sum_{n=0}^{\infty}t(1,1,6;32n+19)q^{2n+1}
=16\cdot 2q(\phq{12}\psq 8+q^2\phq 4\psq{24})\psq 2 =32q\psq 2^2\psq
6.$$ Therefore,
$$\align &\sum_{n=0}^{\infty}t(1,1,6;32n+3)q^n=
16(\phq 2\phq{6}+4q^2\psq 4\psq{12})\psi(q),\tag 3.10
\\&\sum_{n=0}^{\infty}t(1,1,6;32n+19)q^n=32\psi(q)^2\psq 3.
\tag 3.11\endalign$$ From (3.11) and (1.7) we find that
$t(1,1,6;32n+19)=4t(1,1,3;n)$. By (1.6) and (1.10),
$$\sum_{n=0}^{\infty}N(1,1,6;n)q^n
=\varphi(q)^2\phq 6=(\phq 2^2+4q\psq 4^2)\phq 6. $$
Thus,
$$\sum_{n=0}^{\infty}N(1,1,6;2n)q^{2n}=\varphi(q^2)^2\phq 6
=\sum_{n=0}^{\infty}N(1,1,3;n)q^{2n}$$ and so
$N(1,1,6;2n)=N(1,1,3;n)$. On the other hand,
$$\align &\sum_{n=0}^{\infty}N(1,1,3;n)q^n\\&=\varphi(q)^2\phq 3
=(\phq 4^2+4q^2\psq 8^2+4q\psq 4^2)(\phq{12}+2q^3\psq{24})
\\&=\phq 4^2\phq{12}+8q^4\psq 4^2\psq{24}+4q\phq{12}\psq
4^2\\&\q+8q^5\psq 8^2\psq{24}+4q^2\phq{12}\psq 8^2+2q^3\phq
4^2\psq{24}.
\endalign$$
Thus,
$$\align &\sum_{n=0}^{\infty}N(1,1,3;4n)q^{4n}=\phq 4^2\phq{12}+8q^4\psq
4^2\psq{24},
\\&\sum_{n=0}^{\infty}N(1,1,3;4n+1)q^{4n+1}=4q\phq{12}\psq
4^2+8q^5\psq 8^2\psq{24},
\\&\sum_{n=0}^{\infty}N(1,1,3;4n+2)q^{4n+2}=4q^2\phq{12}\psq 8^2,
\\&\sum_{n=0}^{\infty}N(1,1,3;4n+3)q^{4n+3}=2q^3\phq
4^2\psq{24}.
\endalign$$
Substituting $q$ with $q^{1/4}$ yields
$$\align
&\sum_{n=0}^{\infty}N(1,1,3;4n)q^n=\varphi(q)^2\phq{3}+8q\psi(q)^2\psq{6},\tag
3.12
\\&\sum_{n=0}^{\infty}N(1,1,3;4n+1)q^n=4\phq{3}\psi(q)^2+8q\psq
2^2\psq{6},\tag 3.13
\\&\sum_{n=0}^{\infty}N(1,1,3;4n+2)q^n=4\phq{3}\psq 2^2,\tag 3.14
\\&\sum_{n=0}^{\infty}N(1,1,3;4n+3)q^n=2\varphi(q)^2\psq{6}.\tag
3.15
\endalign$$
Combining (3.2) with (3.14) gives $t(1,1,6;4n+1)=4N(1,1,3;4n+2)$,
and combining (3.3) with (3.15) yields
$t(1,1,6;4n+2)=4N(1,1,3;4n+3)$. By (3.13),
$$\align &\sum_{n=0}^{\infty}N(1,1,3;4n+1)q^n
\\&=4\varphi(q)\phq 3\psq 2+8q\psq 2^2\psq 6
\\&=4\psq 2(\phq 4+2q\psq 8)(\phq {12}+2q^3\psq{24})+8q\psq 2^2\psq
6.\endalign$$ Thus,
$$\sum_{n=0}^{\infty}N(1,1,3;8n+1)q^{2n}
=4\psq 2(\phq 4\phq{12}+4q^4\psq 8\psq{24})$$ and
$$\align&\sum_{n=0}^{\infty}N(1,1,3;8n+5)q^{2n+1}\\&=4\psq 2(2q\phq {12}\psq
8+2q^3\phq 4\psq{24}+8q\psq 2^2\psq 6
\\&=4\psq 2\cdot 2q\psq 2\psq 6+8q\psq 2^2\psq 6=16q\psq 2^2\psq 6.
\endalign$$
Therefore,
$$\sum_{n=0}^{\infty}N(1,1,3;8n+1)q^n
=4\psi(q)(\phq 2\phq{6}+4q^2\psq 4\psq{12})\tag 3.16$$ and
$$\sum_{n=0}^{\infty}N(1,1,3;8n+5)q^n=16\psi(q)^2\psq 3.\tag 3.17$$
Combining (3.5) with (3.16) gives $t(1,1,6;8n)=2N(1,1,3;8n+1)$, and
combining (3.6) with (3.17) yields $t(1,1,6;8n+4)=N(1,1,3;8n+5)$. By
(3.12),
$$\align&\sum_{n=0}^{\infty}N(1,1,3;4n)q^n
\\&=\varphi(q)^2\phq 3+8q\varphi(q)\psq 2\psq 6
\\&=(\phq 4^2+4q^2\psq 8^2+4q\psq 4^2)(\phq{12}+2q^3\psq{24})
\\&\qq+8q(\phq 4+2q\psq 8)(\phq{12}\psq 8+q^2\phq 4\psq{24}).
\endalign$$
Thus,
$$\align
&\sum_{n=0}^{\infty}N(1,1,3;4(4n+1))q^{4n+1}=8q^5\psq 8^2\psq{24}
+4q\phq{12}\psq 4^2+8q\phq{12}\phq 4\psq
8\\&\qq\qq\qq\qq\qq\qq\q=8q^5\psq 8^2\psq{24}+12q\phq {12}\psq 4^2,
\\&\sum_{n=0}^{\infty}N(1,1,3;4(4n+2))q^{4n+2}
=4q^2\psq 8^2\phq{12}+16q^2\psq 8^2\phq {12}=20q^2\phq{12}\psq
8^2,
\\&\sum_{n=0}^{\infty}N(1,1,3;4(4n+3))q^{4n+3}
=2q^3\phq 4^2\psq{24}+8q^3\phq 4^2\psq{24}=10q^3\phq 4^2\psq{24}.
\endalign$$
Therefore,
$$\align
&\sum_{n=0}^{\infty}N(1,1,3;16n+4)q^n=8q\psq 2^2\psq{6}+12\phq
{3}\psi(q)^2,\tag 3.18
\\&\sum_{n=0}^{\infty}N(1,1,3;16n+8)q^n
=20\phq{3}\psq 2^2,\tag 3.19
\\&\sum_{n=0}^{\infty}N(1,1,3;16n+12)q^n
=10\varphi(q)^2\psq{6}.\tag 3.20
\endalign$$
Combining (3.8) with (3.19) yields
$t(1,1,6;16n+7)=\f{8}5N(1,1,3;16n+8)$, and combining (3.9) with
(3.20) gives $t(1,1,6;16n+11)=\f{8}5N(1,1,3;16n+8)$.
\par By (1.8)-(1.11),
$$\align\phq 3\psi(q)^2
&=\phq 3\varphi(q)\psq 2 =(\phq{12}+2q^3\psq{24})(\phq 4+2q\psq
8)\psq 2
\\&=\psq 2(\phq 4\phq{12}+4q^4\psq 8\psq{24})
+2q\psq 2(\phq{12}\psq 8+q^2\phq 4\psq{24})
\\&=\phq 4\phq{12}\psq 2+4q^4\psq 2\psq 8\psq{24}+2q\psq 2^2\psq 6.
\endalign$$
Now from (3.18) and the above we deduce that
$$\align
&\sum_{n=0}^{\infty}N(1,1,3;16n+4)q^n \\&=8q\psq 2^2\psq{6}+12\phq
{3}\psi(q)^2
\\&=12\phq 4\phq{12}\psq 2+48q^4\psq 2\psq 8\psq{24}+32q\psq
2^2\psq 6.\endalign$$ Therefore,
$$\align
&\sum_{n=0}^{\infty}N(1,1,3;32n+4)q^{2n}=12\phq 4\phq{12}\psq
2+48q^4\psq 2\psq 8\psq{24},
\\&\sum_{n=0}^{\infty}N(1,1,3;16(2n+1)+4)q^{2n+1}=32q\psq
2^2\psq 6.\endalign$$ Replacing $q$ with $q^{1/2}$ yields
$$\align
&\sum_{n=0}^{\infty}N(1,1,3;32n+4)q^n=12\phq
2\phq{6}\psi(q)+48q^2\psi(q)\psq 4\psq{12},\tag 3.21
\\&\sum_{n=0}^{\infty}N(1,1,3;32n+20)q^n=32\psi(q)^2
\psq 3.\tag 3.22\endalign$$ Combining (3.10) and (3.21) gives
$t(1,1,6;32n+3)=\f 43N(1,1,3;32n+4)$, and combining (3.11) with
(3.22) yields $t(1,1,6;32n+19)=N(1,1,3;32n+20)$.
\par Summarizing the above proves the theorem.

\section*{4. Formulas for $t(1,3,4;n),\ t(1,3,12;n),\ t(1,3,16;n),$\\
$t(1,3,36;n),t(1,3,48;n)$ and $t(3,4,9;n)$}
 \pro{Lemma 4.1} Let
$a,b\in\Bbb Z^+$ with $2\nmid a$. For $n\in\Bbb Z^+$ we have
$$\align &t(a,3a,4b;4n+6a)=2t(a,12a,b;n),\\& t(a,3a,4b;4n+3a)=2t(3a,4a,b;n),
\\&\sum_{n=0}^{\infty}t(a,3a,4b;4n)q^n=8\phq{6a}\psq a\psq b,
\\&\sum_{n=0}^{\infty}t(a,3a,4b;4n+a)q^n=8\phq{2a}\psq{3a}\psq b.
\endalign$$
\endpro
Proof. By (1.9) and (1.11),
$$\align &\sum_{n=0}^{\infty}t(a,3a,4b;n)q^n
=8\psq a\psq{3a}\psq{4b}
\\&=8(\phq{6a}\psq{4a}+q^a\phq{2a}\psq{12a})\psq{4b}
\\&=8(\phq{24a}+2q^{6a}\psq{48a})\psq{4a}\psq{4b}
+q^a(\phq{8a}+2q^{2a}\psq{16a})\psq{12a}\psq{4b}
\\&=8\psq{4b}(\phq{24a}\psq{4a}+q^a\phq{8a}\psq{12a}\\&\qq+2q^{6a}\psq{48a}\psq{4a}
+2q^{3a}\psq{16a}\psq{12a}).\endalign$$ Thus,
$$\align
&\sum_{n=0}^{\infty}t(a,3a,4b;4n)q^{4n}=8\phq{24a}\psq{4a}\psq{4b},
\\&\sum_{n=0}^{\infty}t(a,3a,4b;4n+a)q^{4n+a}=8q^a\phq{8a}\psq{12a}\psq{4b},
\\&\sum_{n=0}^{\infty}t(a,3a,4b;4n+6a)q^{4n+6a}=16q^{6a}\psq{48a}\psq{4a}\psq{4b},
\\&\sum_{n=0}^{\infty}t(a,3a,4b;4n+3a)q^{4n+3a}=16q^{3a}\psq{16a}\psq{12a}\psq{4b}.
\endalign$$
Substituting $q$ with $q^{1/4}$ and then applying (1.7) yields the
result.

\pro{Lemma 4.2} Let $a,b\in\Bbb Z^+$ with $2\nmid a$. If $2\mid b$,
then
$$\align &\sum_{n=0}^{\infty}N(a,3a,4b;8n+5a)q^n=4q^a\psq a\psq{12a}\phq{b/2},
\\&\sum_{n=0}^{\infty}N(a,3a,4b;8n+7a)q^n=4\psq{3a}\psq{4a}\phq{b/2}.
\endalign$$ If $2\nmid b$, then
$$\align &\sum_{n=0}^{\infty}N(a,3a,4b;8n)q^n
\\&\q=\phq{2a}\phq{6a}\phq{2b}+4q^{2a}\psq{4a}\psq{12a}\phq{2b}
\\&\qq+12q^{(a+b)/2}(\phq{6a}\psq{4a}+q^a\phq{2a}\psq{12a})\psq{4b},
\\&\sum_{n=0}^{\infty}N(a,3a,4b;8n+4)q^n
\\&\q=2q^{(b-1)/2}\phq{2a}\phq{6a}\psq{4b}+8q^{2a+(b-1)/2}\psq{4a}\psq{12a}\psq{4b}
\\&\qq+6q^{(a-1)/2}(\phq{6a}\psq{4a}+q^a\phq{2a}\psq{12a})\phq{2b}.
\endalign$$
\endpro
Proof. By (1.13),
$$\align &\sum_{n=0}^{\infty}N(a,3a,4b;n)q^n=\phq a\phq{3a}\phq{4b}
\\&=(\phq {16a}\phq{48a}+4q^{16a}\psq{32a}\psq{96a}+2q^a \phq {48a}\psq
{8a}+2q^{3a}\phq{16a}\psq{24a} \\&\q+6q^{4a}\psq
{8a}\psq{24a}+4q^{13a}\psq
{8a}\psq{96a}+4q^{7a}\psq{24a}\psq{32a})\phq{4b}.\endalign$$ Thus,
$$\align &\sum_{n=0}^{\infty}N(a,3a,4b;4n)q^{4n}
\\&\qq=(\phq {16a}\phq{48a}+4q^{16a}\psq{32a}\psq{96a}+6q^{4a}\psq
{8a}\psq{24a})\phq{4b},
\\&\sum_{n=0}^{\infty}N(a,3a,4b;8n+5a)q^{8n+5a}=4q^{13a}\psq
{8a}\psq{96a}\phq{4b},
\\&\sum_{n=0}^{\infty}N(a,3a,4b;8n+7a)q^{8n+7a}=4q^{7a}\psq{24a}\psq{32a}
\phq{4b}\endalign$$ Substituting $q$ with $q^{1/8}$ in the last two
formulas yields the result in the case $2\mid b$. Now suppose
$2\nmid b$. Substituting $q$ with $q^{1/4}$ in the first formula
gives
$$\align &\sum_{n=0}^{\infty}N(a,3a,4b;4n)q^n
\\&\qq=(\phq {4a}\phq{12a}+4q^{4a}\psq{8a}\psq{24a}+6q^{a}\psq
{2a}\psq{6a})\phq{b} \\&=(\phq
{4a}\phq{12a}+4q^{4a}\psq{8a}\psq{24a}+6q^{a}\psq {2a}\psq{6a})
(\phq{4b}+2q^b\psq{8b}).
\endalign$$
Thus,
$$\align &\sum_{n=0}^{\infty}N(a,3a,4b;8n)q^{2n}
\\&\q=\phq{4a}\phq{12a}\phq{4b}+4q^{4a}\psq{8a}\psq{24a}\phq{4b}
+12q^{a+b}\psq{2a}\psq{6a}\psq{8b},
\\&\sum_{n=0}^{\infty}N(a,3a,4b;8n+4)q^{2n+1}
\\&\q=2q^b\phq{4a}\phq{12a}\psq{8b}+8q^{4a+b}\psq{8a}\psq{24a}\psq{8b}
+6q^a\psq{2a}\psq{6a}\phq{4b}.
\endalign$$
Replacing $q$ with $q^{1/2}$ yields
$$\align &\sum_{n=0}^{\infty}N(a,3a,4b;8n)q^n
\\&=\phq{2a}\phq{6a}\phq{2b}+4q^{2a}\psq{4a}\psq{12a}\phq{2b}
+12q^{(a+b)/2}\psq{a}\psq{3a}\psq{4b},
\\&\sum_{n=0}^{\infty}N(a,3a,4b;8n+4)q^n
\\&=2q^{(b-1)/2}\phq{2a}\phq{6a}\psq{4b}+8q^{2a+(b-1)/2}\psq{4a}\psq{12a}\psq{4b}
\\&\qq+6q^{(a-1)/2}\psq{a}\psq{3a}\phq{2b}.
\endalign$$
By (1.11), $\psq a\psq{3a}=\phq{6a}\psq{4a}+q^a\phq{2a}\psq{12a}$.
Thus the result in the case $2\nmid b$ follows from the above. The
proof is now complete.

\pro{Theorem 4.1} For $n\in\Bbb Z^+$ we have
$$\align &t(1,3,16;n)=2N(1,3,16;2n+5)=t(2,2,3;(n-1)/4)\qtq{for}n\e 1\mod 4,
\\&t(1,3,48;n)=2N(1,3,48;2n+13)=t(1,6,6;n/4)\qtq{for}n\e 0\mod 4,
\\&t(1,3,4;n)=\f 43N(1,3,4;2n+2)\qtq{for}n\e 3,5\mod 8,
\\&t(1,3,12;n)=\f 43N(1,3,12;2n+4)\qtq{for}n\e 0,2\mod 8,
\\&t(1,3,36;n)=\f 43N(1,3,36;2n+10)\qtq{for}n\e 1,7\mod 8,
\\&t(3,4,9;n)=\f 43N(3,4,9;2n+4)\qtq{for}n\e 2,4\mod 8.
\endalign$$
\endpro
Proof. By Lemma 4.1 and (1.8),
$$\sum_{n=0}^{\infty}t(1,3,16;4n+1)q^n=8\phq 2\psq 3\psq 4=8\psq
2^2\psq 3=\sum_{n=0}^{\infty}t(2,2,3;n)q^n.$$ Thus,
$t(1,3,16;4n+1)=t(2,2,3;n)$. By Lemma 4.2,
$$\sum_{n=0}^{\infty}N(1,3,16;8n+7)q^n=4\psq 3\psq 4\phq 2=4\psq
2^2\psq 3.$$ Thus, $t(1,3,16;4n+1)=t(2,2,3;n)=2N(1,3,16;8n+7)$. By
Lemma 4.1,
$$\sum_{n=0}^{\infty}t(1,3,48;4n)q^n=8\phq
6\psi(q)\psq{12}=8\psi(q)\psq
6^2=\sum_{n=0}^{\infty}t(1,6,6;n)q^n.$$ Thus,
$t(1,3,48;4n)=t(1,6,6;n)$.  On the other hand, using Lemma 4.2 we
see that
$$\sum_{n=0}^{\infty}N(1,3,48;8n+5)q^n=4q\psi(q)\psq{12}\phq
6=4q\psi(q)\psq 6^2.$$ Hence
$$\sum_{n=0}^{\infty}N(1,3,48;8n+13)q^n=4\psi(q)\psq 6^2
=\f 12\sum_{n=0}^{\infty}t(1,6,6;n)q^n.$$ Therefore,
$t(1,3,48;4n)=t(1,6,6;n)=2N(1,3,48;8n+13)$.
\par By Lemma 4.1, $t(1,3,4;8n+3)=2t(1,3,4;2n)$.
By (1.11),
$$\sum_{n=0}^{\infty}t(1,3,4;n)q^n=8\psi(q)\psq 3\psq 4
=8(\phq 6\psq 4+q\phq 2\psq{12})\psq 4.$$ Thus,
$$\sum_{n=0}^{\infty}t(1,3,4;2n)q^{2n}=8\phq 6\psq 4^2$$
and so
$$\sum_{n=0}^{\infty}t(1,3,4;8n+3)q^n=2\sum_{n=0}^{\infty}t(1,3,4;2n)q^n
=16\phq 3\psq 2^2.$$ On the other hand, using Lemma 4.2 we see that
$$\align&\sum_{n=0}^{\infty}N(1,3,4;8n)q^n
\\&=\phq 2^2\phq 6+4q^2\phq 2\psq 4\psq{12}+12q\phq 6\psq
4^2+12q^2\phq 2\psq{4}\psq{12}.\endalign$$ Thus,
$$\sum_{n=0}^{\infty}N(1,3,4;16n+8)q^{2n+1}=12q\phq 6\psq 4^2$$
and so
$$\sum_{n=0}^{\infty}N(1,3,4;16n+8)q^n=12\phq 3\psq 2^2=\f {12}{16}
\sum_{n=0}^{\infty}t(1,3,4;8n+3)q^n,$$ which yields
$t(1,3,4;8n+3)=\f 43N(1,3,4;16n+8)$. By Lemma 4.1 and (1.11),
$$\sum_{n=0}^{\infty}t(1,3,4;4n+1)q^n=8\phq 2\psi(q)\psq 3
=8\phq 2(\phq 6\psq 4+q\phq 2\psq{12}).$$ Thus,
$$t(1,3,4;8n+5)q^{2n+1}=8q\phq 2^2\psq{12}.$$ On the other hand,
from Lemma 4.2 we know that
$$\align\sum_{n=0}^{\infty}N(1,3,4;8n+4)q^n&=2\phq 2\phq 6\psq
4+8q^2\psq 4^2\psq{12}
\\&\q+6(\phq 6\psq 4+q\phq 2\psq{12})\phq 2.
\endalign$$
Thus,
$$\sum_{n=0}^{\infty}N(1,3,4;16n+12)q^{2n+1}=6q\phq 2^2\psq{12}
=\f 68\sum_{n=0}^{\infty}t(1,3,4;8n+5)q^{2n+1}$$ and so
$t(1,3,4;8n+5)=\f 43N(1,3,4;16n+12)$.
\par By the above, the formulas for
$t(1,3,16;n),t(1,3,48;n)$ and $t(1,3,4;n)$ are true. Using Lemmas
4.1-4.2 and similar arguments we may deduce the remaining results
for $t(1,3,12;n),$ $t(1,3,36;n)$ and $t(3,4,9;n)$.

\section*{5. Formulas for $t(1,4,7;n),\ t(1,7,12;n),\ t(1,7,28;n),$\\ $t(3,4,21;n)$ and
$t(3,21,28;n)$}
\par We begin with an identity involving $\psi(q)\psq 7$.
\pro{Lemma 5.1} For $q|<1$ we have
$$\align &\psi(q)\psq 7\\&=\psq 8\phq{28}+q^6\phq
4\psq{56}+q\psq{16}\phq{56}+q^3\psq 4\psq{28}+q^{13}\phq 8\psq{112}.
\endalign$$\endpro
Proof. By [Be, p.315],
$$\psi(q)\psq 7=\psq 8\phq{28}+q\psq 2\psq{14}+q^6\phq
4\psq{56}.\tag 5.1$$ Thus,
$$\psq 2\psq {14}=\psq {16}\phq{56}+q^2\psq 4\psq{28}+q^{12}\phq
8\psq{112}.\tag 5.2$$ Combining (5.1) with (5.2) gives the result.

 \pro{Lemma 5.2} Suppose
$a,b\in\Bbb Z^+$ with $2\nmid a$. Then
$$\sum_{n=0}^{\infty}t(a,7a,4b;4n+1)q^n
=8q^{\f{a-1}2}\phq{14a}\psq{4a}\psq
b+8q^{\f{7a-1}2}\phq{2a}\psq{28a}\psq b.$$
\endpro
Proof. Note that $\sum_{n=0}^{\infty}t(a,7a,4b;n)q^n=8\psq
a\psq{7a}\psq{4b}.$ Using Lemma 5.1 we see that
$$\align &\sum_{n=0}^{\infty}t(a,7a,4b;4n+1)q^{4n+1}
\\&=8(q^{a}\phq{56a}\psq{16a} +q^{13a}\phq{8a}\psq{112a})\psq{4b}.
\endalign$$ Replacing $q$ with $q^{1/4}$ yields the result.

 \pro{Theorem 5.1} For $n\in\Bbb Z^+$ we have
$$\align &t(1,4,7;n)=2N(1,4,7;2n+3)\q \t{for}\ n\e 1\mod
4,\\ &t(1,7,12;n)=2N(1,7,12;2n+5)\q\t{for}\ n\e 3\mod 4,
\\&t(1,7,28;n)=2N(1,7,28;2n+9)\q \t{for}\ n\e 1\mod
4,
\\&t(3,4,21;n)=2N(3,4,21;2n+7)\q \t{for}\ n\e 1\mod
4,
\\&t(3,21,28;n)=2N(3,21,28;2n+13)\q \t{for}\ n\e 1\mod
4.
\endalign$$
\endpro
Proof. We only prove the formula for $t(1,4,7;n)$. The other
formulas can be proved similarly. By (1.9),
$$\align &\sum_{n=0}^{\infty}N(1,4,7;n)q^n=\varphi(q)\phq 7\phq 4
\\&=(\phq {16}+2q^4\psq{32}+2q\psq
8)(\phq{112}+2q^{28}\psq{224}+2q^7\psq{56})(\phq{16}+2q^4\psq{32}).
\endalign$$
Thus,
$$\sum_{n=0}^{\infty}N(1,4,7;8n+5)q^{8n+5}
=2q\psq 8(\phq{112}\cdot
2q^4\psq{32}+2q^{28}\psq{224}\cdot\phq{16}). $$
 Replacing $q$ with
$q^{1/8}$ gives
$$\sum_{n=0}^{\infty}N(1,4,7;8n+5)q^n=4\phq{14}\psi(q)\psq
4+4q^3\phq 2\psi(q)\psq{28}.$$ By Lemma 5.2,
$$\sum_{n=0}^{\infty}t(1,4,7;4n+1)q^n=8\phq{14}\psi(q)\psq
4+8q^3\phq 2\psi(q)\psq{28}.$$ Thus
$t(1,4,7;4n+1)=2N(1,4,7;2(4n+1)+3)$. This completes the proof.

\section*{6. Formulas for $t(3,5,4b;n)$ and $t(1,15,4b;n)$ }
\par In order to deduce formulas for $t(3,5,4b;n)$ and
$t(1,15,4b;n)$, we begin with two useful identities.
\pro{Lemma 6.1}
For $|q|<1$ we have
$$\align &\psq 3\psq 5=\phq {60}\psq 8+q^3\psq 2\psq{30}+q^{14}\phq
4\psq{120},\tag 6.1
\\&\psi(q)\psq{15}=\psq
6\psq{10}+q\phq{20}\psq{24}+q^3\phq{12}\psq{40}.\tag 6.2
\endalign$$\endpro
\par We note that (6.1) can be found in [Be,p.377], and (6.2) was
deduced by Xia and Zhang in [XZ].

\pro{Lemma 6.2} For $|q|<1$ we have
$$\align&\psq 3\psq 5=\phq{60}\psq
8+q^{14}\phq 4\psq{120}+q^3\psq{12}\psq{20}
\\&\qq\qq\qq+q^5\phq{40}\psq{48}+q^9\phq{24}\psq{80},
\\&\psi(q)\psq{15}=\phq{120}\psq{16}+q^{28}\phq
8\psq{240}+q^6\psq
4\psq{60}\\&\qq\qq\qq+q\phq{20}\psq{24}+q^3\phq{12}\psq{40}.
\endalign$$
\endpro
Proof. By (6.2), $$\psi(q^2)\psq{30}=\psq
{12}\psq{20}+q^2\phq{40}\psq{48}+q^6\phq{24}\psq{80}.$$ This
together with (6.1) yields the first identity. By (6.1), $$ \psq
6\psq {10}=\phq {120}\psq {16}+q^6\psq 4\psq{60}+q^{28}\phq
8\psq{240}.$$ Combining this with (6.2) gives the second identity.
\pro{Lemma 6.3} Let $b\in\Bbb Z^+$. Then $t(3,5,4b;4n+3)=t(3,5,b;n)$
for $n\in\Bbb Z^+$, and
$$\sum_{n=0}^{\infty}t(3,5,4b;4n+5)q^n=8\psq
b(\phq{10}\psq{12}+q\phq 6\psq{20}).\tag 6.3$$
\endpro
Proof. Note that $\sum_{n=0}^{\infty}t(3,5,4b;n)q^n=8\psq 3\psq
5\psq{4b}$. In view of Lemma 6.2 we see that
$$\sum_{n=0}^{\infty}t(3,5,4b;4n+5)q^{4n+5}=
8\psq{4b}(q^5\phq{40}\psq{48}+q^9\phq{24}\psq{80}).$$ Replacing $q$
with $q^{1/4}$ yields (6.3). Similarly,
$$\sum_{n=0}^{\infty}t(3,5,4b;4n+3)q^{4n+3}=8\psq{4b}\cdot q^3
\psq{12}\psq{20}$$ and so
$$\sum_{n=0}^{\infty}t(3,5,4b;4n+3)q^n=8\psq b\psq 3\psq
5=\sum_{n=0}^{\infty}t(3,5,b;n)q^n,$$ which yields
$t(3,5,4b;4n+3)=t(3,5,b;n)$.  So the lemma is proved.

\pro{Lemma 6.4} Suppose $b\in\Bbb Z^+$ with $2\nmid b$. Then
$$\sum_{n=0}^{\infty}N(3,5,4b;8n+9)q^n=4q^{\f{b-1}2}\phq 6\psq
5\psq{4b}+4q\phq{2b}\psq 5\psq{12}.$$
\endpro
Proof. By (1.9),
$$\align&\sum_{n=0}^{\infty}N(3,5,4b;n)q^n=\phq 3\phq 5\phq {4b}
\\&=(\phq{48}+2q^{12}\psq{96}+2q^3\psq{24})(\phq{80}+2q^{20}\psq{160}+2q^5\psq{40})
\\&\qq\times(\phq{16b}+2q^{4b}\psq{32b}).
\endalign$$
Thus,
$$\align&\sum_{n=0}^{\infty}N(3,5,4b;8n+9)q^{8n+9}
\\&=\phq{48}\cdot 2q^5\psq{40}\cdot 2q^{4b}\psq{32b}
+2q^{12}\psq{96}\cdot 2q^5\psq{40}\cdot \phq{16b}.
\endalign$$
Now replacing $q$ with $q^{1/8}$ yields the result.

 \pro{Theorem 6.1} For $n\in\Bbb Z^+$ we have
$$\align &t(3,4,5;n)=2N(3,4,5;2n+3)\qtq{for}n\e 3\mod 4,
\\&t(3,5,20;n)=2N(3,5,20;2n+7)\qtq{for} n\e 1\mod
4,
\\&t(3,5,36;n)=2N(3,5,36;2n+11)\qtq{for}n\e 3\mod 4,
\\&t(3,5,12;n)=\f 12 N(3,5,12;8n+20)\qtq{for}n\e 1\mod 4.
\endalign$$
\endpro
Proof. By Lemma 6.4,
$$\align &\sum_{n=0}^{\infty}N(3,4,5;8n+9)q^n
\\&=4\psq 5(\phq 6\psq 4+q\phq 2\psq{12})
=4\psi(q)\psq 3\psq 5
\\&=\f 12\sum_{n=0}^{\infty}t(1,3,5;n)q^n.\endalign$$
This together with Lemma 6.3 (with $b=1$) yields
$$2N(3,4,5;8n+9)=t(1,3,5;n)=t(3,4,5;4n+3).\tag 6.4$$
By Lemma 6.3,
$$\sum_{n=0}^{\infty}t(3,5,20;4n+5)q^n=8\psq
5(\phq{10}\psq{12}+q\phq 6\psq{20}).$$ By Lemma 6.4 and the above,
$$\align &\sum_{n=0}^{\infty}N(3,5,20;8n+9)q^n
\\&=4q\psq
5(\phq{10}\psq{12}+q\phq 6\psq{20})=\f
12\sum_{n=0}^{\infty}t(3,5,20;4n+5)q^{n+1}.
\endalign$$
Hence $t(3,5,20;4n+5)=2N(3,5,20;8n+17)$. The remaining results for
$t(3,5,36;n)$ and $t(3,5,12;n)$ can be proved similarly.

\pro{Theorem 6.2} For $n\in\Bbb Z^+$ we have
$$\align&t(1,4,15;n)=2N(1,4,15;2n+5)\qtq{for}n\e 0\mod 4,
\\&t(1,12,15;n)=2N(1,12,15;2n+7)\qtq{for}n\e 2\mod 4,
\\&t(1,15,20;n)=2N(1,15,20;2n+9)\qtq{for}n\e 2\mod 4,
\\&t(3,4,45;n)=2N(3,4,45;2n+13)\qtq{for}n\e 2\mod 4.
\endalign$$
\endpro
Proof. Suppose $a,b\in\Bbb Z^+$ with $2\nmid ab$. Then
$$\sum_{n=0}^{\infty}t(a,15a,4b;n)q^n=8\psq a\psq{15a}\psq{4b}.$$
Applying Lemma 6.2 we see that
$$\align
&\sum_{n=0}^{\infty}t(a,15a,4b;4n)q^{4n}
=8\psq{4b}(\phq{120a}\psq{16a}+q^{28a}\phq{8a}\psq{240a}),
\\&\sum_{n=0}^{\infty}t(a,15a,4b;4n+2)q^{4n+2}
=8\psq{4b}\cdot q^{6a}\psq{4a}\psq{60a}.\endalign$$ Substituting $q$
with $q^{1/4}$ yields
$$\align
&\sum_{n=0}^{\infty}t(a,15a,4b;4n)q^n
=8\psq{b}(\phq{30a}\psq{4a}+q^{7a}\phq{2a}\psq{60a}),\tag 6.5
\\&\sum_{n=0}^{\infty}t(a,15a,4b;4n+2)q^n
=8q^{\f{3a-1}2}\psq{a}\psq{15a}\psq b.\tag 6.6\endalign$$ On the
other hand,
$$\align&\sum_{n=0}^{\infty}N(a,15a,4b;n)q^n=\phq a\phq{15a}\phq{4b}
\\&=(\phq{16a}+2q^{4a}\phq{32a}+2q^a\psq{8a})(\phq{240a}+2q^{60a}\psq{480a}+2q^{15a}
\psq{120a})
\\&\qq\times (\phq{16b}+2q^{4b}\psq{32b}).
\endalign$$
Thus,
$$\align &\sum_{n=0}^{\infty}N(a,15a,4b;8n+3a)q^{8n+3a}
\\&\q=2q^{15a}\psq{120a}(2q^{4a}\psq{32a}\phq{16b}+2q^{4b}\psq{32b}\phq{16a}),
\\&\sum_{n=0}^{\infty}N(a,15a,4b;8n+a+4)q^{8n+a+4}
\\&\q=2q^a\psq{8a}(2q^{60a}\psq{480a}\phq{16b}+2q^{4b}\psq{32b}\phq{240a}).
\endalign$$
Replacing $q$ with $q^{1/8}$ yields
$$\aligned
&\sum_{n=0}^{\infty}N(a,15a,4b;8n+3a)q^n=4q^{2a}\psq{4a}\psq{15a}\phq{2b}
+4q^{\f{3a+b}2}\psq{15a}\psq{4b}\phq{2a},
\\&\sum_{n=0}^{\infty}N(a,15a,4b;8n+a+4)q^n=4q^{\f{15a-1}2}\psq a\psq{60a}\phq{2b}
+4q^{\f{b-1}2}\psq a\psq{4b}\phq{30a}.
\endaligned\tag 6.7$$
Taking $a=b=1$ in (6.5) and (6.7) gives
$$\align \sum_{n=0}^{\infty}t(1,4,15;4n)q^n&=8\psi(q)(\phq{30}\psq
4+q^7\phq
2\psq{60})\\&=2\sum_{n=0}^{\infty}N(1,4,15;8n+5)q^n,\endalign$$
which implies $t(1,4,15;4n)=2N(1,4,15;8n+5)$. Taking $a=1$ and $b=3$
in (6.6) yields
$$\sum_{n=0}^{\infty}t(1,12,15;4n+2)q^n=8q\psi(q)\psq{15}\psq 3.$$
Taking $a=1$ and $b=3$ in (6.7) yields
$$\sum_{n=0}^{\infty}N(1,12,15;8n+3)q^n=4q^2\psq{15}(\phq 6\psq
4+q\phq 2\psq{12})=4q^2\psi(q)\psq 3\psq{15}.$$ Thus,
$t(1,12,15;4n+2)=2N(1,12,15;8n+11)$. The remaining parts for
$t(1,15,20;n)$ and $t(3,4,45;n)$ can be proved similarly.

\section*{7. Some conjectures on $t(a,b,c;n)$}
\par Let $a,b,c,n\in\Bbb Z^+$. Based on calculations on Maple, in this section
we pose many conjectures on $t(a,b,c;n)$.

\pro{Conjecture 7.1} Let $n\in\Bbb Z^+$ with $n\not\e
47,63\mod{64}$. Then
$$t(1,2,5;n)=\cases 4N(1,2,5;2n+2)&\t{if $n\e 0,1\mod 4$,}
\\N(1,2,5;2n+2)&\t{if $n\e 2\mod 8$ or $n\e 11\mod{32}$,}
\\2N(1,2,5;2n+2)&\t{if $n\e 6\mod 8$,}
\\\f 8{5}N(1,2,5;2n+2)&\t{if $n\e 3,7\mod {16}$,}
\\\f 43N(1,2,5;2n+2)&\t{if $n\e 27\mod {32}$,}
\\\f{16}{13}N(1,2,5;2n+2)&\t{if $n\e 15,31\mod {64}$.}
\endcases$$
\endpro

\pro{Conjecture 7.2} Let $n\in\Bbb Z^+$ with $n\not\e 14\mod{16}$.
Then
$$t(1,5,10;n)
=\cases 4N(1,5,10;2n+4)&\t{if $n\e 0,3\mod 4$,}\\
2N(1,5,10;2n+4)&\t{if $n\e 1\mod 8$,}
\\ N(1,5,10;2n+4)&\t{if $n\e 5\mod 8$ or $n\e 26\mod{32}$,}
\\\f 8{5}N(1,5,10;2n+4)&\t{if $n\e 2,6\mod {16}$,}
\\\f 43N(1,5,10;2n+4)&\t{if $n\e 10\mod {32}$,}
\endcases$$
\endpro

\pro{Conjecture 7.3} Let $n\in\Bbb Z^+$ with $n\not\e 14\mod{16}$.
Then
 $$t(1,6,9;n)=\cases N(1,6,9;2n+4)&\t{if $n\e 3\mod 8$ or $n\e 18\mod {32}$,}
\\2N(1,6,9;2n+4)&\t{if $n\e 7\mod 8$,}
\\4N(1,6,9;2n+4)&\t{if $n\e 0,1\mod 4$,}
\\\f 8{5}N(1,6,9;2n+4)&\t{if $n\e 6,10\mod {16}$,}
\\\f 43N(1,6,9;2n+4)&\t{if $n\e 2\mod{32}$.}
\endcases$$
 \endpro

 \pro{Conjecture 7.4} Let $n\in\Bbb Z^+$ and
$$r(n)=t(1,4,4;n)-\f 12N(1,4,4;8n+9).$$
\par $(\t{\rm i})$ If $n\e 1\mod 3$, then $r(n)=0$.
\par $(\t{\rm ii})$ If $n\e 0\mod 3$, then
$$r(n)=\cases (-1)^{k-1}(6k-3)&\t{if $n=\f{9k^2-9k}2$ for $k\in\Bbb
Z$,}\\0&\t{otherwise.}\endcases$$
\par $(\t{\rm iii})$ If $n\e -1\mod 3$, then
$$r(n)=\cases (-1)^k(6k-1)&\t{if $n+1=\f{9k^2-3k}2$ for $k\in\Bbb
Z$,}
\\(-1)^{k-1}(6k+1)&\t{if $n+1=\f{9k^2+3k}2$ for $k\in\Bbb
Z$,}
\\0&\t{otherwise.}\endcases$$
\endpro

\pro{Conjecture 7.5} Let $n\in\Bbb Z^+$ with $n\e 1\mod 8$. Then
$$\align&t(1,1,30;n)=2N(1,1,30;2n+8),
\\&t(2,3,11;n)=2N(2,3,11;2n+4),
\\&t(2,7,7;n)=2N(2,7,7;2n+4),
\\&t(3,5,12;n)=2N(3,5,12;2n+5).\endalign$$
\endpro

\pro{Conjecture 7.6} Let $n\in\Bbb Z^+$ with $n\e 5\mod 8$. Then
$$\align &t(1,2,13;n)=2N(1,2,13;2n+4),
\\&t(1,14,49;n)=2N(1,14,49;2n+16),
\\&t(2,5,9;n)=2N(2,5,9;2n+4),
\\&t(3,3,10;n)=2N(3,3,10;2n+4).
\endalign$$
\endpro

\pro{Conjecture 7.7} Let $n\in\Bbb Z^+$ with $n\e 3\mod 8$. Then
$$\align &t(1,1,14;n)=2N(1,1,14;2n+4),
\\&t(1,7,8;n)=\f 13N(1,7,8;8n+16),
\\&t(1,21,42;n)=2N(1,21,42;2n+16),
\\&t(2,5,25;n)=2N(2,5,25;2n+8),
\\&t(2,9,21;n)=2N(2,9,21;2n+8),
\\&t(3,5,8;n)=\f 13N(3,5,8;8n+16),
\\&t(3,6,7;n)=2N(3,6,7;2n+4).
\endalign$$
\endpro

\pro{Conjecture 7.8} Let $n\in\Bbb Z^+$ with $n\e 7\mod 8$. Then
$$\align
&t(1,18,45;n)=2N(1,18,45;2n+16),
\\&t(2,3,27;n)=2N(2,3,27;2n+8),
\\&t(2,15,15;n)=2N(2,15,15;2n+8),
\\&t(5,5,6;n)=2N(5,5,6;2n+4).\endalign$$
\endpro

\pro{Conjecture 7.9} Let $n\in\Bbb Z^+$ with $n\e 0\mod 8$. Then
$$\align
 &t(1,5,50;n)=2N(1,5,50;2n+14),
\\&t(1,9,30;n)=2N(1,9,30;2n+10).\endalign$$
\endpro

\pro{Conjecture 7.10} Let $n\in\Bbb Z^+$ with $n\e 4\mod 8$. Then
$$\align
&t(1,2,21;n)=2N(1,2,21;2n+6),
\\&t(1,6,33;n)=2N(1,6,33;2n+10),
\\&t(1,10,45;n)=2N(1,10,45;2n+14),
\\&t(2,11,11;n)=N(2,11,11;2n+6).
\endalign$$
\endpro

\pro{Conjecture 7.11} Let $n\in\Bbb Z^+$ with $n\e 2\mod 8$. Then
$$\align &t(1,1,22;n)=N(1,1,22;2n+6),
\\&t(1,15,24;n)=\f 13N(1,15,24;8n+40),
\\&t(1,15,40;n)=\f 13N(1,15,40;8n+56),
\\&t(1,18,21;n)=2N(1,18,21;2n+10),
\\&t(1,22,33;n)=2N(1,22,33;2n+14),
\\&t(5,6,45;n)=2N(5,6,45;2n+14).\endalign$$
\endpro

\pro{Conjecture 7.12} Let $n\in\Bbb Z^+$ with $n\e 6\mod 8$. Then
$$\align &t(1,7,32;n)=\f 12N(1,7,32;8n+40),
\\&t(1,13,26;n)=2N(1,13,26;2n+10),
\\&t(1,25,30;n)=2N(1,25,30;2n+14),
\\&t(3,5,32;n)=\f 12N(3,5,32;8n+40),
\\&t(3,7,14;n)=2N(3,7,14;2n+6)
.\endalign
$$
\endpro


\end{document}